\documentclass[aos,MSNbibl,nameyear,dvips]{arximspdf}
\usepackage{graphicx}

\doi{10.1214/12-AOS979} 
\referstodoi{10.1214/11-AOS949}
\volume{40}
\issue{4}
\pubyear{2012}
\firstpage{1968}
\lastpage{1972}

\makeatletter
\newproclaim{efx}{Expectation step (E step)}
\newproclaim{mfa}{Maximization step (M step)}
\def\bE{\mathbb{E}}
\def\bR{\mathbb{R}}
\def\rank{\operatorname{rank}}
\def\tr{\operatorname{trace}}
\def\var{\operatorname{Var}}
\def\argmin{\mathop{\operatorname{arg\,min}}}
\newlength{\dhatheight}
\newcommand{\dhat}[1]{%
\settoheight{\dhatheight}{\ensuremath{\hat{#1}}}%
\addtolength{\dhatheight}{-0.25ex}%
\hspace*{1pt}\hat{\hspace*{-1pt}\vphantom{\rule{1pt}{\dhatheight}}%
\smash{\hat{#1}}}}
\makeatother

\begin{document}
\begin{frontmatter}

\title{Discussion: Latent variable graphical model selection via
convex optimization\thanksref{T1}}

\thankstext{T1}{Supported in part by NSF Career Award DMS-08-46234.}
\runtitle{Comment}

\begin{aug}
\author[A]{\fnms{Ming} \snm{Yuan}\corref{}\ead[label=e1]{myuan@isye.gatech.edu}}

\runauthor{M. Yuan}
\affiliation{Georgia Institute of Technology}
\address[A]{H. Milton Stewart School of Industrial\\
\quad and Systems Engineering\\
Georgia Institute of Technology\\
Atlanta, Georgia 30332\\
USA} 
\end{aug}

\received{\smonth{2} \syear{2012}}



\end{frontmatter}

I want to start by congratulating Professors Chandrasekaran, Parrilo
and Willsky for this fine piece of work. Their paper, hereafter
referred to as CPW,
addresses
one of the biggest practical challenges of Gaussian graphical
models---how to make inferences for a graphical model in the presence
of missing variables. The difficulty comes from the fact that the
validity of conditional independence relationships implied by a
graphical model relies critically on the assumption that all
conditional variables are observed, which of course can be unrealistic.
As CPW shows, this is not as hopeless as it might appear to be. They
characterize conditions under which a conditional graphical model can
be identified, and offer a penalized likelihood method to reconstruct
it. CPW notes that with missing variables, the concentration matrix of
the observables can be expressed as the difference between a sparse
matrix and a low-rank matrix; and suggests to exploit the sparsity
using an $\ell_1$ penalty and the low-rank structure by a trace norm
penalty. In particular, the trace norm penalty or, more generally,
nuclear norm penalties, can be viewed as a convex relaxation to the
more direct rank constraint. Its use oftentimes comes as a necessity
because rank constrained optimization could be computationally
prohibitive. Interestingly, as I note here, the current problem
actually lends itself to efficient algorithms in dealing with the rank
constraint, and therefore allows for an attractive alternative to the
approach of CPW.






\section{Rank constrained latent variable graphical Lasso}

Recall that the penalized likelihood estimate of CPW is defined as
\[
(\hat{S}_n,\hat{L}_n)=\argmin_{L\succeq0, S-L\succ0}
\bigl\{-\ell(S-L,\Sigma_O^n)+\lambda_n\bigl(\gamma\|S\|_{1}+\tr(L)\bigr)\bigr\},
\]
where the vector $\ell_1$ norm and trace/nuclear norm penalties are
designated to induce sparsity among\vadjust{\goodbreak} elements of $S$ and low-rank
structure of $L$ respectively. Of course, we can attempt a more direct
rank penalty as opposed to the nuclear norm penalty on $L$, leading to
\[
(\hat{S}_n,\hat{L}_n)=\argmin_{L\succeq0, S-L\succ0}
\bigl\{-\ell(S-L,\Sigma_O^n)+\lambda_n\bigl(\gamma\|S\|_{1}+\rank(L)\bigr)\bigr\};
\]
or for computational purposes, it is more convenient to consider the
constrained version:
\[
(\dhat{S}_n,\dhat{L}_n)=\mathop{\argmin_{L\succeq0, S-L\succ
0}}_{\rank(L)\le r} \{-\ell(S-L,\Sigma_O^n)+\lambda_n\|S^\dag\|
_{1}\},
\]
for some integer $0\le r\le p$, where $S^\dag=S-{\rm diag}(S)$, that
is, $S^\dag$ equals $S$ except that its diagonals are replaced by $0$.
This slight modification reflects our intention to encourage sparsity
on the off-diagonal entries of $S$ only. The remaining discussion,
however, can be easily adapted to deal with the original vector
$\ell_1$ penalty on $S$. It is clear that when $r=0$, that is, $L=0$,
this new estimator reduces to the so-called graphical Lasso estimate
(\verb+glasso+, for short) of Yuan and Lin (\citeyear{YuaLin07}). See
also Banerjee, El Ghaoui and d'Aspremont (\citeyear{BanElGdAs08}),
Friedman, Hastie and Tibshirani (\citeyear{FriHasTib08}), and Rothman
et al. (\citeyear{Rotetal08}). Drawn to this similarity, I shall
hereafter refer to this method as the latent variable graphical Lasso,
or \verb+LVglasso+, for short.

Common wisdom on $(\dhat{S}_n,\dhat{L}_n)$ is that it is infeasible to
compute because of the nonconvexity of the rank constraint.
Interestingly, though, this more direct approach actually allows for
fast computation, thanks to a combination of EM algorithm and some
recent advances in computing graphical Lasso estimates for
high-dimensional problems.

\section{An EM algorithm}

The constraint $\rank(L)\le r$ amounts to postulating $r$ latent
variables. The latent variable model naturally has a missing data
formulation. It is clear that when observing the complete data
$X=(X_O^{\sf T},X_H^{\sf T})^{\sf T}$, the \verb+LVglasso+ estimator
becomes
\[
\hat{K}_\lambda=\argmin_{K\in\bR^{(p+r)\times(p+r)}, K\succ0}
\{L(K)+\lambda\|K_{O}^\dag\|_{1}\},
\]
where
\[
L(K)=-\ln\det(K)+\tr\bigl(\Sigma_{(O H)}^nK\bigr)
\]
and $\Sigma_{(O H)}^n$ is the sample covariance matrix of the full
data. Now that $X_H$ is unobservable, we can use an EM algorithm which
iteratively applies the following two steps:

\begin{efx*}
Calculate the expected value of the
penalized negative log-likelihood function, with respect to the
conditional distribution of $X_H$ given $X_O$\vadjust{\goodbreak} under the current
estimate $K^{(t)}$ of $K$, leading to the so-called Q function:
\begin{eqnarray*}
Q\bigl(K|K^{(t)}\bigr)&=&\bE_{X_H|X_O,K^{(t)}}[L(K)+\lambda\|K_{O}^\dag\|
_{1}]\\
&=&-\ln\det(K)+\tr\bigl\{\bE_{X_H|X_O,K^{(t)}}\bigl(\Sigma_{(O
H)}^n\bigr)K\bigr\}+\lambda\|K_{O}^\dag\|_{1}.
\end{eqnarray*}
Recall that $X_H|X_O,K^{(t)}$ follows a normal distribution with
\[
\bE\bigl(X_H|X_O,K^{(t)}\bigr)=\Sigma^{(t)}_{HO}\bigl(\Sigma^{(t)}_{O}\bigr)^{-1}X_O
\]
and
\[
\var\bigl(X_H|X_O,K^{(t)}\bigr)=\Sigma^{(t)}_{H}-\Sigma^{(t)}_{HO}\bigl(\Sigma
^{(t)}_{O}\bigr)^{-1}\Sigma^{(t)}_{OH},
\]
where $\Sigma^{(t)}=(K^{(t)})^{-1}$. Therefore,
\[
\bE_{X_H|X_O,K^{(t)}}(\Sigma_{OH}^n)=\Sigma^n_{O}\bigl(\Sigma
^{(t)}_{O}\bigr)^{-1}\Sigma^{(t)}_{oH}
\]
and
\[
\bE_{X_H|X_O,K^{(t)}}(\Sigma_{H}^n)=\Sigma^{(t)}_{H}-\Sigma
^{(t)}_{HO}\bigl(\Sigma^{(t)}_{O}\bigr)^{-1}\Sigma^{(t)}_{OH}+\Sigma
^{(t)}_{HO}\bigl(\Sigma^{(t)}_{O}\bigr)^{-1}\Sigma^n_{O}\bigl(\Sigma
^{(t)}_{O}\bigr)^{-1}\Sigma^{(t)}_{OH}.
\]
\end{efx*}

\begin{mfa*}
Maximize $Q(\cdot|K^{(t)})$ over all
$(p+r)\times(p+r)$ positive definite matrices. We first note that if we
replace the penalty term $\|K_O^\dag\|_1$ with $\|K^\dag\| _1$, then
maximizing $Q(\cdot|K^{(t)})$ becomes a \verb+glasso+ problem:
\[
\max_{K\in\bR^{(p+r)\times(p+r)}, K\succ0}
\bigl\{-\ln\det(K)+\tr\{WK\}+\lambda\|K^\dag\|_{1}\bigr\},
\]
where $W=\bE_{X_H|X_O,K^{(t)}}(\Sigma_{(O H)}^n)$. As shown in
Banerjee, El~Ghaoui and\break
d'Aspremont (\citeyear{BanElGdAs08}),
\citeauthor{FriHasTib08} (\citeyear{FriHasTib08}) and Yuan
(\citeyear{Yua08}), this problem can be solved iteratively. At each
iteration, one row and, correspondingly, one column of $K$, due to
symmetry, are updated by solving a Lasso problem. The same idea can be
applied here to maximize $Q(\cdot |K^{(t)})$. The only difference is
that in each of the Lasso problems, we leave the coordinates
corresponding to the latent variables unpenalized. This extension has
\begin{figure}
\includegraphics{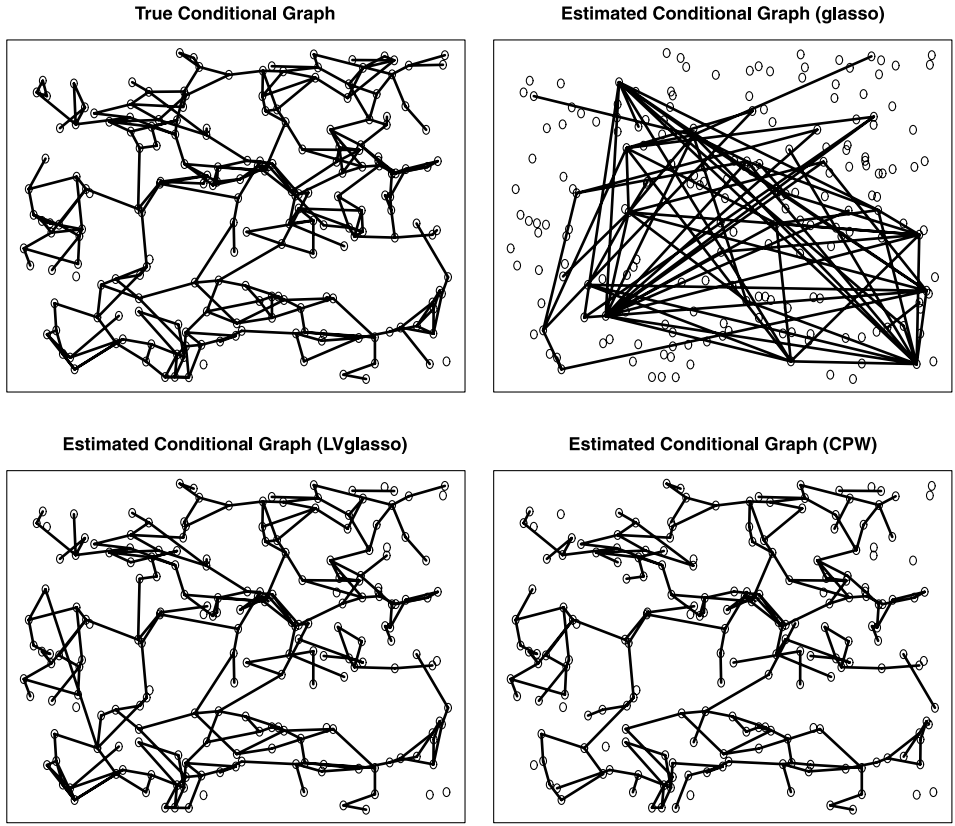}
\caption{True graphical model and its estimates.}\label{fig:graphs}
\end{figure}
been implemented in the \verb+R+ package \verb+glasso+
[\citeauthor{FriHasTib08} (\citeyear{FriHasTib08})].
\end{mfa*}

\section{Example}

For illustration purposes, I conducted a simple numerical experiment.
In this experiment the interest was in recovering a $p=198$ dimensional
graphical model with $h=2$ missing variables. The graphical model was
generated in a similar fashion as that from Meinshausen and B\"uhlmann
(\citeyear{MeiBuh06}). I first simulated 198 locations uniformly over a
square.
\begin{figure}
\includegraphics{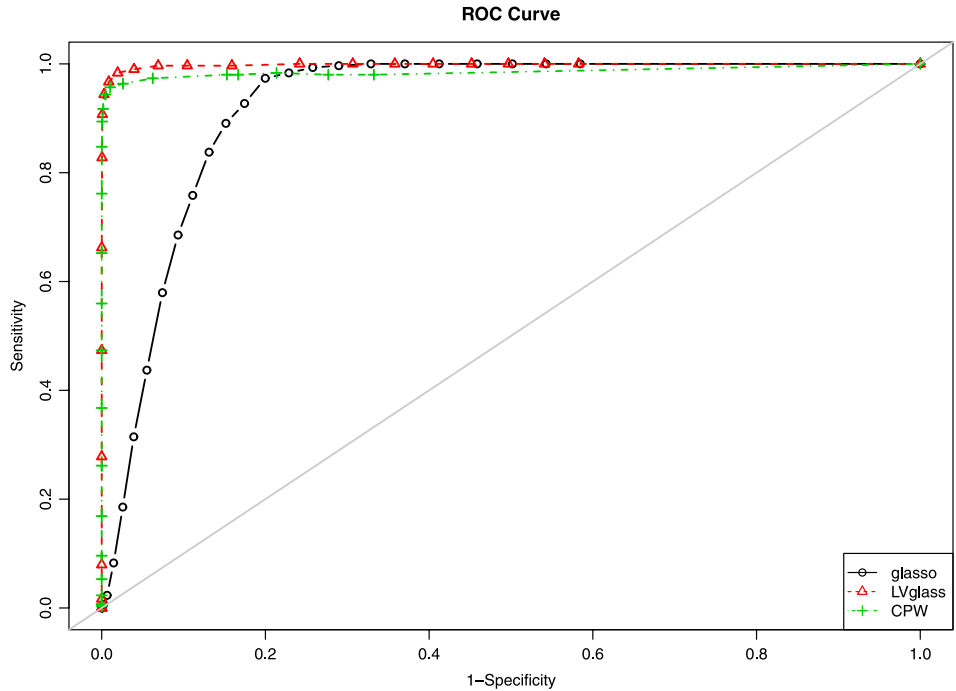}
\caption{Accuracy of reconstructed conditional graphical model.}\label{fig:roc}
\end{figure}
Between each pair of locations, I put an edge with probability
$2\phi(d\sqrt{p})$, where $\phi(\cdot)$ is the density function of the
standard normal distribution and $d$ is the distance between the two
locations, unless one of the locations is already connected with four
other locations. The two hidden variables were connected with all $p$
observed variables. The entries of the inverse covariance matrix
corresponding to the edges between the observables were assigned with
value $0.2$, between the observables and the latent variables were
assigned with a uniform random value between 0 and $0.12$, to ensure
the positive definiteness. A typical simulated graphical model among
the 198 observed variables conditional on the two latent variables is
given in the top left panel of Figure~\ref{fig:graphs}. We apply both
the method of CPW and \verb+LVglasso+, along with \verb+glasso+, to the
data. We used the \verb+MATLAB+ code provided by CPW to compute their
estimates. As observed by CPW, their estimate typically is insensitive
to a wide range of values of $\gamma$, and we report here the results
with the default choice of $\gamma=5$ without loss of generality.
Similarly, for \verb+LVglasso+, little variation was observed for
$r=2,\ldots, 10$, and we shall focus on $r=2$ for brevity. The choice
of $\lambda$ plays a critical role for both methods. We compute both
estimators for a fine grid of $\lambda$. With the main focus on
recovering the conditional graphical model, that is, the sparsity
pattern of~$S$, we report in Figure~\ref{fig:roc} the ROC curve for
both methods. For contrast, we also reported the result for
\verb+glasso+ which neglects the missingness. In Figure~\ref{fig:graphs}, we also presented the estimated graphical model for
each method that is closest to the truth. These results clearly
demonstrate the necessity of accounting for the latent variables. It is
also interesting to note that the rank constrained estimator performs
slightly better in this example over the trace norm penalization method
of CPW.


The preliminary results presented here suggest that direct rank
constraint may provide a competitive alternative to the trace norm
penalization for recovering graphical models with latent variables. It
is of interest to investigate more rigorously how the two methods
compare with each other.


%

\printaddresses

\end{document}